\documentclass[fleqn]{mat01}
\usepackage{times,mathtimy,amssymb,boxit}
\begin{document}
\setcounter{page}{281}
\firstpage{281}

\newcommand{\n}{\noindent}
\newcommand{\bb}[1]{\mathbb{#1}}
\newcommand{\cl}[1]{\mathcal{#1}}
\newcommand{\inner}[2]{\mbox{$\Big \langle #1,#2\Big \rangle$}}
\newcommand{\kf}[1]{{K(\cdot,#1)}}%
\newcommand{\kzf}[1]{{K_0(\cdot,#1)}}%
\def\x{\mbox{$\pmb{x}$}}%

\def\defi{\trivlist\item[\hskip\labelsep{DEFINITION.}]}
\def\remark{\trivlist\item[\hskip\labelsep{\it Remark.}]}
\def\remaa{\trivlist\item[\hskip\labelsep{\it {\rm [}Remark.}]}
\def\remarks{\trivlist\item[\hskip\labelsep{\it Remarks}]}
\def\noot{\trivlist\item[\hskip\labelsep{{\it Note.}}]}
\def\leem{\trivlist\item[\hskip\labelsep{{\it Lemma.}}]}
\def\hyp{\trivlist\item[\hskip\labelsep{\it Hypothesis.}]}
\def\nnthm{\trivlist\item[\hskip\labelsep{{\bf Theorem.}}]}

\font\xx=msam5 at 10pt
\def\ab{\mbox{\xx{\char'03}}}

\title{Equivalence of quotient Hilbert modules}

\markboth{Ronald G Douglas and Gadadhar Misra}{Equivalence of quotient Hilbert modules}

\author{RONALD G DOUGLAS and GADADHAR MISRA$^{*}$}

\address{Texas A\&M University, College Station, Texas 77843, USA\\
\noindent $^{*}$Indian Statistical Institute, R.V. College Post, Bangalore~560~059, India\\
\noindent E-mail: rgd@tamu.edu; gm@isibang.ac.in}


\volume{113}

\mon{August}

\parts{3}

\Date{MS received 15 June 2002}

\begin{abstract}
Let $\cl{M}$ be a Hilbert module of holomorphic functions over a
natural function algebra $\mathcal{A}(\Omega)$, where $\Omega
\subseteq \bb{C}^m$ is a bounded domain.  Let $\cl{M}_0\subseteq
\cl{M}$ be the submodule of functions vanishing to order $k$ on a
hypersurface $\cl{Z} \subseteq \Omega$.  We describe a method, which
in principle may be used, to construct a set of complete unitary
invariants for quotient modules $\cl{Q}=\cl{M} \ominus \cl{M}_0$. The
invariants are given explicitly in the particular case of $k = 2$.
\end{abstract}


\maketitle

\section{Preliminaries}

Let $\Omega$ be a bounded domain in $\bb{C}^m$ and $\cl{Z}\subseteq
\Omega$ be an analytic hypersurface defined (at least, locally) as the
zero set of a single analytic function $\varphi$.  Let
$\cl{A}(\Omega)$ be the algebra of functions obtained by taking the
closure with respect to the supremum norm on $\Omega$ of all functions
which are holomorphic on a neighbourhood of $\Omega$.  Let $\cl{M}$ be
a Hilbert space consisting of holomorphic functions on $\Omega$.  We
assume that the evaluation functionals $h \to h(w), ~h\in \cl{M},
~w\in \Omega$ are bounded.  This ensures, via the Riesz representation
theorem, that there is a unique vector $K(\cdot, w)\in \cl{M}$
satisfying the reproducing property
\begin{equation*}
h(w) = \langle h, K(\cdot, w) \rangle, ~h \in \cl{M},~w\in \Omega.
\end{equation*}
In this paper, a module $\cl{M}$ over the
function algebra $\cl{A}(\Omega)$ will consist of a Hilbert space
$\cl{M}$ as above together with a continuous action of the algebra
$\cl{A}(\Omega)$ in the sense of (\cite{D-P}, Definition 1.2).
Suppose, we are given a quotient
module $\cl{Q}$ over the function algebra $\cl{A}(\Omega)$. This amounts
to the existence of a resolution of the form
\begin{equation} \label{resolution}
0 \longleftarrow \cl{Q} \longleftarrow \cl{M} \longleftarrow
\cl{M}_0 \longleftarrow 0,
\end{equation}
where $\cl{M}_0 \subseteq \cl{M}$ are
both modules over the algebra $\cl{A}(\Omega)$. We make the additional
assumption that the submodule $\cl{M}_0$ consists of functions in
$\cl{M}$ which vanish to some fixed order $k$ on the hypersurface
$\cl{Z}$.  Then (cf. \cite{D-M-V}, (1.5)) the module $\cl{M}_0$ may be
described as
\begin{equation*}
\cl{M}_0 = \left\lbrace f\in \cl{M}: \frac{\partial^{\ell}f}{\partial
z_1^{\ell}} (z) = 0,\:\: z\in U\cap \cl{Z},\:\: 0\leq
\ell \leq k-1 \right\rbrace,
\end{equation*}
where $U$ is some open subset of $\Omega$.

Let $\partial$ denote the differentiation along the unit normal to the
hypersurface $\cl{Z}$.  Recall (cf. \cite{D-M-V}) that the map $J: \cl{M}
\to \cl{M} \otimes \bb{C}^k$ defined by
\begin{equation*}
h \mapsto (h, \partial h, \partial^2 h, \ldots , \partial^{k-1} h), ~
h\in \cl{M}
\end{equation*}
plays a crucial role in identifying the quotient module. The
requirement that
\begin{equation*}
\{(e_n, \partial e_n, \ldots , \partial^{k-1}e_n)_{n \geq 0}:
(e_n)_{n\geq 0} \mbox {~is an orthonormal basis in $\cl{M}$~}\}
\end{equation*} is
an orthonormal basis in $\mbox{ran}\: J$, makes the map $J$ unitary
onto its range $J\cl{M} \subseteq \cl{M}\otimes \bb{C}^k$.  Thus we
obtain a pair of modules $J\cl{M}_0$ and $J\cl{M}$, where $J\cl{M}_0$
is the submodule of all functions in $J\cl{M}$ which vanish on
$\cl{Z}$.  In this realisation, the module $J\cl{M}$ consists of
holomorphic functions taking values in $\bb{C}^k$.
Let $\bb{C}^{k\times k}$ denote the linear space of all $k\times k$
matrices over the field of complex numbers.
We recall that a function
$K:\Omega\times \Omega \to \bb{C}^{k\times k}$ satisfying
\begin{equation} \label{existence reprod}
\sum_{i,j=1}^n \inner{K(\omega_i,\omega_j)\zeta_j}{\zeta_i}_E ~\geq~0 ,~~w_1,
\ldots,\omega_n\in \Omega, ~~\zeta_1,\ldots,\zeta_n \in E, n\geq 0
\end{equation}
is said to be a {\em nonnegative definite $($nnd$)$ kernel} on $\Omega$.
Given such an nnd kernel $K$ on $\Omega$, it is easy to construct a
Hilbert space
$\mathcal{M}$ of functions on $\Omega$ taking values in
$\bb{C}^{k\times k}$ with the
property
\begin{equation} \label{reproducing property}
\inner{f(\omega)}{\zeta}_{\bb{C}^k} = \inner{f}{\kf{\omega}\zeta},~ w\in
\Omega,~\zeta\in \bb{C}^k,~f\in \mathcal{M}.
\end{equation}
The Hilbert space $\mathcal{M}$ is simply the completion of the linear
span of all vectors of the form $\kf{\omega}\zeta$, $\omega\in
\Omega$, $\zeta\in \bb{C}^k$, with inner product defined by
(\ref{reproducing property}).  Conversely, let $\mathcal{M}$ be a
Hilbert space of functions on $\Omega$ taking values in $\bb{C}^k$.  Let
$e_\omega : \mathcal{M} \to \bb{C}^k$ be the evaluation functional defined
by $e_\omega(f) = f(\omega)$, $\omega\in \Omega$, $f\in
\mathcal{M}$.  If $e_\omega$ is bounded for each $\omega \in \Omega$,
then it is easy to verify that the Hilbert space $\mathcal{M}$
possesses a reproducing kernel
$K(z,\omega)= e_z e_\omega^*$, that is, $K(z,\omega)\zeta \in
\mathcal{M}$ for each $\omega\in \Omega$ and $K$ has the reproducing
property (\ref{reproducing property}). Finally,  the reproducing
property (\ref{reproducing property}) determines the reproducing
kernel $K$ uniquely.  If $e_n$ is an orthonormal basis in
$\mathcal{M}$ then it is not hard to verify that the reproducing
kernel $K$ has the representation
\begin{equation*}
K(z,w) = \sum_{n=0}^\infty e_n(z) e_n(w)^*, ~~z,w \in \Omega,
\end{equation*}
where $e_n(z)$ is thought of as a linear map from $\bb{C}$ to $\bb{C}^k$.
Of course, this sum is independent of the choice of the orthonormal
basis $e_n$ since $K$ is uniquely determined.

The module $J\mathcal{M}$ possesses a reproducing kernel $JK$ in the
sense described above.
It is natural to construct this kernel by forming the sum:
\begin{equation*}
JK(z,w) = \sum_{n=0}^\infty (Je_n)(z) (Je_n)(w)^*,~z,w \in \Omega.
\end{equation*}
This prescription then allows the identification of the reproducing
kernel $JK:\Omega \times \Omega \to \bb{C}^{k\times k}$ for the module $J\cl{M}$:
\begin{equation}\label{jet kernel}
(JK)_{\ell,j}(z,w) =\big (\partial^{\ell}\bar{\partial}^{j}K \big )
(z,w),~~~ 0\leq \ell,j \leq k-1.
\end{equation}
It is then easy to verify, using the unitarity of the map $J$, that
$JK$ has the reproducing property:
\begin{equation*}
\langle h , JK(\cdot , w) \zeta \rangle = \langle h(w) , \zeta
\rangle, ~w\in \Omega,~ \zeta \in \bb{C}^k.
\end{equation*}
The module action for $J\cl{M}$ is defined in a natural manner.  Indeed,
let $Jf$ be the array
\begin{equation}\label{jet action}
(Jf)_{\ell,j} = \begin{cases}
\begin{pmatrix}\ell\\  j\end{pmatrix} (\partial^{\ell-j}f), & 0\leq \ell\leq j \leq k-1\\
0, & \mbox{otherwise} \end{cases}
\end{equation}
for $f \in \cl{A}(\Omega)$.
We may now define the module action to be $J_f: h \to Jf \cdot Jh$.
Notice that $Jf$ is a $k\times k$ matrix-valued function on $\Omega$
while $J_f$ is the module action, that is, it is an operator on
$J\mathcal{M}$.
The action of the adjoint is then easily seen to be
\begin{equation} \label{adjoint action}
J_f^* J\kf{w} \cdot {\bf x} = J\kf{w} (Jf)(w)^*\cdot {\bf x}, ~{\bf x} \in \bb{C}^k.
\end{equation}

We will say that two modules over the algebra $\cl{A}(\Omega)$ are {\em
isomorphic} if there exists a unitary module map between them.

It is shown in \cite{D-M-V} that the quotient module $\cl{Q}$ is
isomorphic to  $J\cl{M}\ominus J\cl{M}_0$.  Once this is done, we are reduced
to the multiplicity free case. Thus our previous results from \cite{D-M2}
apply and we conclude that the quotient module $\cl{Q}$ is
the restriction of $J\cl{M}$ to the hypersurface $\cl{Z}$.

Let $\cl{M}$ be any Hilbert module over the function algebra
$\cl{A}(\Omega)$.  In particular, each of the coordinate functions
$z_i, ~1\leq i \leq m$ in $\bb{C}^m$ acts boundedly as the multiplication
operator $M_i$ on $\cl{M}$. Let {\bf M} denote this commuting $m$-tuple of
multiplication operators.  We denote by ${\bf M}^*$ the $m$-tuple
$(M_1^*,\ldots ,M_m^*)$.  To each $m$-tuple {\bf M}, we
associate the operator $D_{\bf M}:\mathcal{M} \to \cl{M}\otimes
\bb{C}^k$ defined by $D_{\bf M}h = (M_1h,\ldots ,M_mh)$,
$h\in\cl{M}$.

The class $B_n(\Omega)$ was introduced in \cite{C-D} for
a single operator.  This definition was then adapted to the general case of
an $m$-tuple of commuting operators (cf. \cite{C-S}).
We let $\Omega^*\subseteq \bb{C}^m$ denote the domain $\{w\in
\bb{C}^m: \bar{w} \in \Omega\}$ and say that ${\bf M}^*$ is in
$B_k(\Omega^*)$ if

\begin{enumerate}
\renewcommand{\labelenumi}{(\roman{enumi})}
\leftskip .5pc
\item $\mbox{Ran~}D_{{\bf M}^* - w}$ is closed for all
$w\in \Omega^*$,
\item span~$\{\ker D_{{\bf M}^* - w}: w\in \Omega^*\}$
is dense in $\cl{M}$,
\item $\dim~\ker D_{{\bf M}^* - w}=n$
for all $w\in \Omega^*$,
\end{enumerate}
where ${\bf M}^* - w = (M_1^* - w_1, \ldots ,M^*_m-w_m)$.

If the adjoint of the $m$-tuple of multiplication operators is in
$B_n(\Omega^*)$ (for some $n\in \bb{N}$), then we say that $\cl{M}$ is
in $B_n(\Omega^*)$.  The assumption that $\cl{M}$ is in
$B_1(\Omega^*)$ includes, among other things, (a) the existence of a
common eigenvector $\gamma(w) \in \cl{M}$, that is, $M_i^* \gamma(w) =
\bar{w}_i \gamma(w)$, for $w\in \Omega^*$, (b) the dimension of the
common eigenspace at $\bar{w}$ is $1$.  Furthermore, it is possible to
choose $\gamma(w)$ so as to ensure that the map $w \to \gamma(w)$ is
anti-holomorphic.  Thus we obtain an anti-holomorphic hermitian line
bundle $E$ over $\Omega$ whose fiber at $w$ is the one-dimensional
subspace of $\cl{M}$ spanned by the vector $\gamma(w)$, that is,
$\gamma$ is an anti-holomorphic frame for $E$.  In the case of $n > 1$,
a similar construction of an anti-holomorphic hermitian vector bundle of
rank $n$ can be given.  In our case, it is easy to verify that
$K(\cdot, w)$, the reproducing kernel at $w$, is a common eigenvector
for the $m$-tuple $(M_1^*, \ldots ,M_m^*)$.  Since $K(\cdot,w)$ is
anti-holomorphic in the second variable, it provides a natural frame
for the associated bundle $E$.  The metric with respect to this frame
is obviously the real analytic function $K(w,w)$.

Before we continue, we make the additional assumption that the module
$\cl{M}$, which occurs in the resolution (\ref{resolution}) of the
quotient module $\cl{Q}$, lies in the class $B_1(\Omega^*)$.  Let
$i:\cl{Z} \to \Omega$ be the inclusion map and $i^*:\cl{A}(\Omega) \to
\cl{A}(\cl{Z})$ be the pullback.  Then $\cl{Q}$ is clearly also a
module over the smaller algebra $i^*\big (\cl{A}(\Omega) \big )$.  We
identify this latter algebra with $\cl{A}(\cl{Z})$.  Let $(\cl{Q},
\cl{A}(\cl{Z}))$ stand for $\cl{Q}$ thought of as a module over the
smaller algebra $\cl{A}(\cl{Z})$.  Although it is possible that
$(\cl{Q}, \cl{A}(\cl{Z}))$ lies in $B_k(\cl{Z}^*)$ whenever $\cl{M}$
is in $B_1(\Omega^*)$, we were able to prove it only in some special
cases (\cite{D-M-V}, Proposition 3.6).  However, in this paper, we
assume that the quotient module $(\cl{Q}, \cl{A}(\cl{Z}))$ always lies
in $B_k(\cl{Z}^*)$.  These assumptions make it possible to associate
(a) an anti-holomorphic hermitian line bundle $E$ over the domain
$\Omega$ with the module $\cl{M}$ and (b) an anti-holomorphic jet
bundle $JE_{|\rm res~\cl{Z}}$ of rank $k$ over the domain $\cl{Z}$
with the module $(\cl{Q}, \cl{A}(\cl{Z}))$.  The details of the {\em
jet} construction are given in (\cite{D-M-V}, pp.~375--377).  One of
the main results in \cite{C-D} states that two modules $\cl{M}$ and
$\tilde{\cl{M}}$ in $B_k(\Omega)$ are isomorphic if and only if the
associated bundles are locally equivalent.  While the local
equivalence of bundles is completely captured in the case of line
bundles by the curvature, it is more complicated in the general case
(cf. \cite{C-D}).  We recall that the quotient module $\cl{Q}$ may be
described completely by specifying the action of the algebra
$\cl{A}_k(\cl{Z}):= \cl{A}(\cl{Z}) \otimes \bb{C}^{k\times k}$
(cf. \cite{D-M-V}, p.~385).  The action of the algebra
$\cl{A}_k(\cl{Z})$, in particular, includes the multiplication induced
by the local defining function $\varphi$, namely,
\begin{equation*}(J\varphi)_{| {\rm res~}\cl{Z}} :J\cl{M}_{| {\rm
res~}\cl{Z}} \to J\cl{M}_{| {\rm res~}\cl{Z}}.
\end{equation*}
To exploit methods of \cite{C-D}, it is better to work with the
adjoint action.  To describe the adjoint action, we first
construct a natural anti-holomorphic frame ({\em not
necessarily} orthonormal) for the jet bundle $E$ on $\Omega$.  Let
$\{\varepsilon_\ell: ~1\leq \ell \leq k\}$ be the standard orthonormal
basis in $\bb{C}^k$.  For a fixed $w\in \Omega$, let
$e_1=\sum_{\ell=1}^k \partial^{\ell-1} K(z, w)\otimes
\varepsilon_\ell$ be
simply the image of $K(z,w)$ in $J\cl{M}$.  It is then clear that
$\{e_j(w): 1\leq j \leq k \}$, where $e_j(w): = (\bar{\partial}^{j-1} e_1)(w)$
is a natural anti-holomorphic frame for $JE$.  (Of course, as is to be
expected, $e_\ell(w)$, $1\leq \ell \leq k$ are the columns of the
reproducing kernel $JK$ given in (\ref{jet kernel}).) Thus the fiber
of the jet bundle $JE$ at $w\in \Omega$ is spanned by the set of
vectors $\{e_\ell(w)\in J\cl{M} : 1\leq \ell \leq k\}$.

Suppose we start with a resolution of the form (\ref{resolution}).
Then we have at our disposal the domain $\Omega\subseteq \bb{C}^m$ and
the hypersurface $\cl{Z} \subseteq \Omega$.  Let $\varphi$ be a local
defining function for $\cl{Z}$ (cf. \cite{D-M-V}, p.~367).  Then
$\varphi$ lies in $\cl{A} (\cl{Z})$ and induces a nilpotent action on
each fiber of the jet bundle $JE_{|\rm res~\cl{Z}}$ via
the map $J_\varphi^*$, that is,
\begin{equation} \label{nilaction}
 (J_\varphi^* e_\ell)(w) = J\kf{w} (J\varphi)(w)^* \varepsilon_\ell.
\end{equation}
Therefore in this picture, with the assumptions we have made along
the way, we see that the {\em quotient modules} $\cl{Q}$ must meet the
requirement listed in (i)--(iii) of the following Definition.

\begin{defi}\label{Gen CD}$\left.\right.$\vspace{.5pc}

\noindent We will say that the module $\mathcal{Q}$ over the algebra
$\cl{A}(\Omega)$ is a {\em quotient module} in the class $B_k(\Omega,
\cl{Z})$ if

\begin{enumerate}
\renewcommand{\labelenumi}{(\roman{enumi})}
\leftskip .5pc
\item there exists a resolution of the module
$\mathcal{Q}$ as in eq.~(\ref{resolution}), where
the module $\cl{M}$ appearing in the resolution is required to be in
$B_1(\Omega^*)$,

\item the module action on $\mathcal{Q}$ translates to
the nilpotent action $J_\varphi$ on $J\cl{M}_{| {\rm res~}\cl{Z}}$
which is an isomorphic copy of $\cl{Q}$,
\item the module $\big (\cl{Q}, \cl{A}(\cl{Z}) \big )$ is
in $B_k(\cl{Z}^*)$.
\end{enumerate}
\end{defi}

In this paper, we obtain a complete set of unitary invariants for a
module $\mathcal{Q}$ in the class $B_2(\Omega, \mathcal{Z})$.  This
means that the module $\mathcal{Q}$ admits a resolution of the form
(\ref{resolution}) and the module $\mathcal{M}$ that appears in this
resolution lies in $B_1(\Omega)$.  However, it is possible to
considerably weaken this latter hypothesis as explained in the Remark
below.

\begin{remaa}
Although we have assumed the module $\mathcal{M}$ to be in the class
$B_1(\Omega)$, it is interesting to note that the proof of our Theorem
requires much less.  Specifically, the requirement that the
`$\mbox{Ran~}D_{{\bf M}^* - w}$ is closed' is necessary to associate
an anti-holomorphic vector bundle with the module.  However, in our
case, there is already a natural anti-holomorphic vector bundle which
is deteremined by the frame $w \to K(\cdot ,w)$.  Indeed, if we assume
that the module $\mathcal{M}$ contains the linear space $\mathcal{P}$
of all the polynomials and $\mathcal{P}$ is dense in $\mathcal{M}$,
then the eigenspace at $w$ is forced to be one dimensional.  (To prove
this, merely note that for any eigenvector $x$ at $w$ and all
polynomials $p$, we have
\begin{equation*}
\langle p , x \rangle = \langle M_p 1, x \rangle
= \langle 1 , M_p^* x \rangle
= p(w) \langle 1 , x \rangle
= \langle p , c K(\cdot , w) \rangle,
\end{equation*}
where $c= \overline{\langle 1 , x \rangle}$.  It follows that $x= c
K(\cdot,w)$.)  Finally, the linear span of the set of eigenvectors
$\{K(\cdot ,w): w \in \Omega\}$ is a dense subspace of the module
$\mathcal{M}$.  Therefore, for our purposes, it is enough to merely
assume that
\begin{enumerate}
\renewcommand{\labelenumi}{(\alph{enumi})}
\item {\em $\mathcal{M}$ is a Hilbert module consisting of
holomorphic functions on $\Omega$,}
\item {\em the module $\mathcal{M}$ contains the linear space of all
polynomials $\mathcal{P}$ and that $\mathcal{P}$ is dense,}
\item {\em $\mathcal{M}$ possesses a reproducing
kernel $K$.}
\end{enumerate}\vspace{-.3pc}

It is then clear that the same holds for the quotient module
$\mathcal{Q}$, where $\mathcal{P}$ consists of $\bb{C}^k$-valued
polynomials and $K$ takes values in $\bb{C}^{k\times k}$. Hence, if $x$
is an eigenvector at $w$ for the module $(\mathcal{Q},
\mathcal{A}(\mathcal{Z}))$, we claim that it belongs to the range of
$K(.,w)$ which is the $k$-dimensional subspace $\{K(\cdot , w) v \in
\mathcal{Q}: v \in \bb{C}^k \}$ of $\mathcal{Q}$.  As before, for
$1\leq j \leq k$, let $\varepsilon_j$ be the standard unit vector in
$\bb{C}^k$ and $p= \sum_{j=1}^k p_j \otimes \varepsilon_j$ be a
$\bb{C}^k$-valued polynomial. Then we have
$\langle p, x \rangle =
\sum_{j=1}^k \langle M_{p_j} \varepsilon_j , x \rangle = \sum_{j=1}^k
\langle \varepsilon_j , M_{p_j}^* x \rangle = \sum_{j=1}^k p_j(w)
\langle \varepsilon_j , x \rangle\\ = \sum_{j=1}^k \langle p, K(\cdot,w)
\varepsilon_j \rangle \langle \varepsilon_j , x \rangle = \langle p,
\sum_{j=1}^k c_j K(.,w) \varepsilon_j \rangle,$
where $c_j = \overline {\langle \varepsilon_j , x \rangle}$.  Thus $x$ is in the
range of $K(\cdot, w)$ as claimed.  Therefore the dimension of
the eigenspace at $w$ equals the dimension of range $K(.,w)$ which is~$k$.]
\end{remaa}\vspace{.5pc}

We now raise the issue of adapting the techniques of \cite{C-D} to
find a complete set of unitary invariants for characterizing the
quotient modules $\cl{Q}$ in the class $B_k(\Omega, \cl{Z})$.  While
the methods described below will certainly yield results in the
general case, we have chosen to give the details of our results in the
case of $k=2$.  The reason for this choice is dictated by the simple
nature of these invariants in this case. Furthermore, these are
extracted out of the curvature and the canonical metric for the bundle
$E$.

\section{Canonical metric and curvature}

Let $\cl{M}$ be a module in $B_1(\Omega^*)$ and the reproducing kernel
$K(\cdot, w)$ be the anti-holomorphic frame for the associated bundle $E$.
If $\tilde{\cl{M}}$ is another module in the class $B_1(\Omega^*)$ with
reproducing kernel $\tilde{K}(\cdot, w)$, then it is clear that any
isomorphism between these modules must map $K(\cdot,w)$ to a multiple
$\psi(w)$ of $\tilde{K}(\cdot, w)$, where $\psi(w)$ is a non
zero complex number for $w\in \Omega$.  Moreover, the map $w \to
\psi(w)$ has to be anti-holomorphic.  It follows that
$\cl{M}$ and $\tilde{\cl{M}}$ are isomorphic if and only if $\tilde{K}(z,w) =
\overline{\psi(z)} K(z,w) \psi(w)$  (cf. \cite{C-S}, Lemma~3.9)
for some anti-holomorphic function $\psi$.
There are two ways in which this ambiguity may be eliminated.

The first approach is to note that if the two modules
$\tilde{\cl{M}}$ and $\cl{M}$ are isomorphic, then $\tilde{K}(z,z)/K(z,z) =
|\psi(z)|^2$.  Since $\psi$ is holomorphic, it follows that
\begin{equation} \label{curvature}
\sum_{i,j=1}^m\partial_i\bar{\partial}_j \log
\big (K(z,z)/\tilde{K}(z,z) \big ) dz_i\wedge d\bar{z}_j= 0.
\end{equation}
On the other hand, if we have two modules for which equation
(\ref{curvature}) holds, then the preceding argument shows that they must
be isomorphic.  It is then possible to find, in a small
simply connected neighbourhood of some fixed point $w_0$, a harmonic
conjugate $v(w)$ of the harmonic function $u(w):= \log
\tilde{K}(w,w)/K(w,w)$.  The new kernel defined by
\ $\tilde{\tilde{\!K}}(z,w) = \exp(u(z)+iv(z))\tilde{K}(z,w)
\overline{\exp(u(w)+iv(w))}$ determines a module
\ \,$\tilde{\tilde{\!\!\cl{M}}}$ isomorphic to $\tilde{M}$ but with the
additional property that the metric \ $\tilde{\tilde{\!K}}(w,w) = K(w,w)$.
It is then easy to see that the map taking $K(\cdot, w)$ to
\ $\tilde{\tilde{\!K}}(\cdot, w)$ extends linearly to an isometric module
map.  Therefore, $\sum_{i,j=1}^m\partial_i\bar{\partial}_j \log
K(z,z) dz_i\wedge d\bar{z}_j$ is a complete invariant for the\break
module $\mathcal{M}$

The second approach is to normalise the reproducing kernel $K$\!, that is,
define the kernel $K_0(z,w) = \psi(z) K(z,w)
\overline{\psi(w)}$, where $\psi(z) =
K(z,w_0)^{-1}K(w_0,w_0)^{1/2}$ for $z$ in some open subset $\Omega_0
\subseteq \Omega$ and some fixed but arbitrary $w_0 \in \Omega_0$.
Also, $\Omega_0$ can be chosen so as to ensure $\psi_{|{\rm
res}~\Omega_0} \not = 0$.
This reproducing kernel determines a module isomorphic to $\cl{M}$ but
with the added property that $K_0(z,w_0)$ is
the constant function $1$.  If $\cl{M}$ and $\tilde{\cl{M}}$ are two
modules in $B_1(\Omega^*)$, then it is shown in (\cite{C-S},
Theorem~4.12) that they are isomorphic if and only if the normalisations
$K_0$ and $\tilde{K}_0$ of the respective reproducing kernels at some
fixed point are equal.  As before, it is then easy to see that the map
taking $K(\cdot, w)$ to \ $\tilde{\tilde{\!K}}(\cdot, w)$ extends linearly
to an isometric module map.  The normalised kernel $K_0$ is therefore
a complete unitary invariant for the module $\mathcal{M}$.

Notice that if a module $\mathcal{M}$ is isomorphic to
$\tilde{\mathcal{M}}$, then the module map $\Gamma$ is induced by a nonvanishing
function $\Phi$ on $\Omega$, that is, $\Gamma=M_\Phi$ (\cite{C-S}, Lemma 3.9).
Consequently, if $\mathcal{M}_0$ is the submodule of functions
vanishing to order $k$ on $\mathcal{Z}$, then $\Gamma(\mathcal{M}_0)$
is the submodule of functions vanishing to order $k$ in
$\tilde{\mathcal{M}}$.  It follows that if $\mathcal{M}$ and
$\tilde{\mathcal{M}}$ are isomorphic modules, then the corresponding
quotient modules must be isomorphic as well.
Therefore we can make the following assumption without any loss of generality.

\begin{hyp}
Now we make a standing hypothesis that the kernel for the module
$\mathcal{M}$ appearing in the resolution of the quotient module
$\mathcal{Q}$ is normalised.
\end{hyp}\vspace{.5pc}

Recall that if $E$ is a hermitian holomorphic vector bundle of rank
$k$ over the domain $\Omega\subseteq \bb{C}^m$, then it is possible to
find a holomorphic frame ${\bf s} = (s_1,\ldots ,s_k)$ such that (a)
$\langle s_i(w_0) , s_j(w_0) \rangle = 1$, (b) $\partial_j \langle
{\bf s}(w), {\bf s}(w) \rangle _{| w = w_0} = 0$ for $1\leq j \leq m$
(cf. \cite{W}, Lemma~2.3).  We offer below a variation of this Lemma
for the jet bundle $JE$ corresponding to the hypersurface
$\mathcal{Z}\subseteq \Omega$ and the Hilbert module $\mathcal{M}$ in
the class $B_1(\Omega)$.  We state the following Lemma in terms of
a frame for the bundle associated with the module $\mathcal{M}$.  There
is an obvious choice for such a frame in terms of the reproducing
kernel of the module.  The relationship
between the reproducing kernel of the module and the
hermitian metric of the associated bundle was explained in
(\cite{D-M-V}, \S~2). Let $\langle {\bf s}(w) , {\bf s}(w_0)
\rangle$ be the matrix of inner products, that is, $\langle {\bf s}(w),
{\bf s}(w_0) \rangle_{ij} = \langle {s}_i(w) , {s}_j(w_0)
\rangle_\mathcal{M}$, $1\leq i,j \leq k$ for some fixed but arbitrary
$w_0\in \mathcal{Z}$ and all $w\in \mathcal{Z}$.

\begin{leem} \label{normal}
{\it Let $\mathcal{M}$ be Hilbert module in $B_1(\Omega)$ and
$\mathcal{M}_0 \subseteq \mathcal{M}$ be the submodule consisting of
functions vanishing on the hypersurface $\mathcal{Z}\subseteq
\Omega$.    Then there exists an anti-holomorphic frame ${\bf s}$ for
the jet bundle JE satisfying
\begin{equation*}
\langle {\bf s}(w) , {\bf s}(w_0) \rangle_{|{\rm res~}\mathcal{Z}} =
\begin{pmatrix} 1 & 0 \\ 0 & S(w) \end{pmatrix},
\end{equation*}
for $w\in \mathcal{Z}$ and some anti-holomorphic function $S$ on $\mathcal{Z}$.}
\end{leem}

\begin{proof}
Let us assume, without loss of generality, that $w_0=0$.  We first
observe that if we replace the module $\mathcal{M}$ by an isomorphic
copy, then the class of the associated bundle {\it JE} does not change.  Indeed, if
$\mathcal{M}$ and $\tilde{\mathcal{M}}$ are isomorphic modules, then
there is an anti-holomorphic map $\varphi$ which induces a metric preserving
bundle map of the associated bundles $E$ and $\tilde{E}$.  It is
then clear that the map $J_\varphi^*$ induces a bundle map of the
corresponding jet bundles.  Therefore, we may assume that the
reproducing kernel $K$ for the module $\mathcal{M}$ is normalised,
that is, $K(z,0)=1$.  Let $(\tilde{z},\tilde{w})$ denote
(temporarily) the normal coordinates in $\Omega\times \Omega$.  From the
expansion
\begin{equation*}
K(z,w) = \sum_{\ell,n=0}^\infty K_{\ell,n}(z,w) \tilde{z}^\ell
\bar{\tilde{w}}^n,\quad z, w\in \mathcal{Z}
\end{equation*}
it is clear that $K_{\ell n}(z,0)=0$ for $\ell\not=0$ and $n=0$.  Since
$K(z,w)= \overline{K(w,z)}$, it follows that $K_{\ell n}(0,w)=0$ for $\ell
=0$ and $n\not =0$.  However, $K_{\ell n}(z,w) = (\partial^\ell
\bar{\partial}^n K)_{|\tilde{z}=0,\tilde{w}=0}(z,w)$.
Hence $(\!(K_{\ell n}(z,w))\!)_{\ell,n=0}^{k-1}=JK_{|{\rm
res~}\mathcal{Z}}(z,w)$ for $z,w \in \mathcal{Z}$ by definition
(\ref{jet kernel}).  Recall that $e_\ell(w) = \sum_{j =1}^k
\bar{\partial}^{\ell-1}\partial^{j -1} K (\cdot ,w) \otimes
\varepsilon_\ell$, for $1\leq \ell \leq k$ is an anti-holomorphic frame
for the jet bundle $JE$.  It follows that $\langle e_\ell(w), e_n(0)
\rangle = (JK)_{\ell n} (0 , w)$.  But $(JK)_{\ell n} (0,w) = K_{\ell n} (0,w) = 0
$ for $\ell =0$ as long as $n\not = 0$.  The proof is completed by taking
${\bf s}(w) = \{e_1(w), \ldots ,e_k (w)\}$.\hfill \ab
\end{proof}

There is a canonical connection $D$ on the bundle $JE$ which is
compatible with the metric and has the property $D^{\prime \prime} =
\bar{\partial}$.  Let $C^\infty_{1,1}(\Omega, E)$ be the space of
$C^\infty$ sections of the bundle $\wedge^{(1,1)}T^*\Omega\otimes E$.
The curvature tensor $\cl{K}$ associated with the canonical connection
$D$ is in $C^\infty_{1,1}(\Omega, \mbox{herm}(E,E))$.  Moreover, if
$h$ is a local representation of the metric in some open set, then
$i\cl{K} = \bar{\partial}({h}^{-1} \partial{h})$.  The holomorphic
tangent bundle $T\Omega_{|{\rm res~ } \cl{Z}}$ naturally splits as
$T\cl{Z} \dot{+} N\cl{Z}$, where $N\cl{Z}$ is the normal bundle and is
realised as the quotient $T\Omega_{|{\rm res~ } \cl{Z}}/T\cl{Z}$.  The
co-normal bundle $N^*\cl{Z}$ is the dual of $N\cl{Z}$; it is the
sub-bundle of $T\Omega_{|{\rm res~ } \cl{Z}}$ consisting of cotangent
vectors that vanish on $T\cl{Z} \subseteq T\Omega_{|{\rm res ~}
\cl{Z}}$.  Indeed, the class of the conormal bundle $N^*\cl{Z}$
coincides with $[-\cl{Z}]_{|{\rm res ~} \cl{Z}}$ via the adjunction
formula I (\cite{G-H}, p.~146).  Let $P_1$ be the projection onto
$N^*\cl{Z}$ and $P_2=(1-P_1)$ be the projection onto $T^*\cl{Z}$.
Now, we have a splitting of the $(1,1)$ forms as follows:
\begin{equation*}
\wedge^{(1,1)}T^* \Omega_{|{\rm res ~}\cl{Z}} = \sum_{i,j=1}^2 P_i\big
(\wedge^{(1,0)}T^*\Omega _{|{\rm res ~}\cl{Z}}\big ) \wedge P_j \big
(\wedge^{(0,1)}T^*\Omega_{|{\rm res ~}\cl{Z}} \big ).
\end{equation*}
Accordingly, we have the component of the curvature along the transversal
direction to $\cl{Z}$ which we denote by $\cl{K}_{\rm trans}$.
Clearly, $\cl{K}_{\rm trans}= (P_1 \otimes I) \cl{K}_{|{\rm res~ }\cl{Z}}$.
Similarly, let the component of the curvature along tangential
directions to $\cl{Z}$ be $\cl{K}_{\rm tan}$.  Again, $ \cl{K}_{\rm
tan} = (P_2 \otimes I)\cl{K}_{|{\rm res~} \cl{Z}}$.
(Here $I$ is the identity map on the vector space $\mbox{herm}(E,E)$.)

Recall that the fiber of the jet bundle $JE_{|{\rm res~}\mathcal{Z}}$ at $w\in
\mathcal{Z}$ is spanned by the set of vectors
$\bar{\partial}^{\ell-1}K(\cdot , w)$, $1\leq \ell \leq k$.  Thus the
module action $J_\varphi^*$ can be determined by calculating it on the
set $\{\bar{\partial}^{\ell-1}K(\cdot , w): 1\leq \ell \leq k \mbox{ and }
w\in \mathcal{Z} \}$.  This calculation is given in eq.~(\ref{nilaction}).  We therefore obtain an anti-holomorphic bundle map
$J_\varphi^*$ on the bundle $JE_{|{\rm res~}\mathcal{Z}}$.  Thus the
isomorphism of two quotient modules in $B_k(\Omega, \mathcal{Z})$
translates to a question of  equivalence of the pair
$(JE_{|{\rm res~}\mathcal{Z}}, J_\varphi^*)$.
This merely amounts to finding an anti-holomorphic bundle map
$\theta:JE_{|{\rm res~}\mathcal{Z}}\to JE_{|{\rm res~}\mathcal{Z}}$
which intertwines $J_\varphi^*$.
It is clear that if we could find such a bundle map $\theta$, then the
line sub-bundles corresponding to the frame $K(\cdot ,w)$, $w\in
\mathcal{Z}$ must be equivalent.  From this it is evident that
the curvatures $\mathcal{K}_{\rm tan}$ in the tangential directions
must be equal.  Also, we can calculate the matrix representation
for the nilpotent action at $w$, as given in (\ref{nilaction}),  with respect to the orthonormal basis
obtained via the Gram--Schmidt process applied to the holomorphic frame
at $w$.  A computation shows that the matrix entries involve the
curvatures $\mathcal{K}_{\rm trans}$ in the transverse direction and
its derivatives.  It is not clear if the intertwining condition can be
stated precisely in terms of these matrix entries.
In the following section we show, as a result of
some explicit calculation, that if $k=2$ then the curvature in the
transverse direction must also be equal.  We also find that an
additional condition must be imposed to determine the
isomorphism class of the quotient\break modules.

\section{The case of rank 2 bundles}

In this case, the adjoint action of $\varphi$ on $\cl{Q} \cong
\left .J\cl{M} \right |_{\rm res~  \cl{Z}}$ produces a nilpotent
bundle map on $JE$ which, at $w\in \cl{Z}$, is described easily:
\begin{equation*}
e(w) := \left (\begin{smallmatrix} K(\cdot, w) \cr {\partial}
K(\cdot, w) \end{smallmatrix} \right ) \to 0 \mbox{ and }
(\bar{\partial}e)(w) := \left (\begin{smallmatrix}
\bar{\partial} K(\cdot, w) \cr
\partial \bar{\partial} K(\cdot, w) \end{smallmatrix} \right ) \to
\overline{(\partial \varphi)(w)} {\rm e}(w)
\end{equation*}
on the spanning set $\{ e(w), (\bar{\partial}e)(w): w\in \cl{Z} \}$
for the fiber $J\!E(w)$ of the jet bundle $J\!E$ at $w\in \cl{Z}$.
Thus the adjoint action induced by $\varphi$ determines a nilpotent
$N(w)$ of
order $2$ defined by $\left (\begin{matrix} 0 &
\overline{(\partial \varphi)(w)} \cr 0 & 0 \end{matrix} \right )$ on each
fiber $J\!E(w)$, $w\in \cl{Z}$ with respect to
the basis $\{ e(w), ({\partial} e)(w) \}$.  Now,
consider the orthonormal basis:
$\{\gamma_0(w), \gamma_1(w)\},$ where
\begin{align*}
\gamma_0(w) &= \|e(w)\|^{-1} e(w),\\
\gamma_1(w) &= a(w) e(w) + b(w) (\bar{\partial}e)(w),~w\in \cl{Z}.
\end{align*}
The coefficients $a(w)\mbox{ and }b(w)$ can be easily calculated
(cf. \cite{C-D}, p.~195):
\begin{align*}
-a(w) \|e(w)\|^3 &= \langle (\partial e)(w), e(w) \rangle
(-\cl{K}_{\rm trans}(w))^{-1/2},\\
b(w) \|e(w)\| &= (- \cl{K}_{\rm trans}(w))^{-1/2},
\end{align*}
where $\cl{K}_{\rm trans}(w)$ denotes the curvature in the transversal
direction.  In the case of a line bundle, we have the following
explicit formula:
\begin{equation} \label{linecurvature}
\cl{K}_{\rm trans}(w) = P_1 \Big (\sum_{i,j=1}^m\partial_i\bar{\partial}_j
\log\|e(w)\|^2 dz_i\wedge d\bar{z}_j \Big ), ~w\in \cl{Z}.
\end{equation}
The nilpotent action $N_{\rm orth}(w)$ at the fiber $J\!E(w), ~w\in
\cl{Z}$ with respect to the
orthonormal basis  $\{\gamma_0(w), \gamma_1(w)\}$ is given by
\begin{equation*}
\left (\begin{matrix} 0 & b(w) \|e(w)\| (\partial \varphi)(w) \cr
0 & 0 \end{matrix} \right ).
\end{equation*}

Now, we are ready to prove the main theorem which gives a complete set
of invariants for quotient modules in the class $B_2(\Omega,
\cl{Z})$. At first, it may appear that the condition {angle} of the
theorem stated below depends on the choice of the holomorphic frame.
But we remind the reader that the normalisation of the kernel
$K$ for the module $\mathcal{M}$ ensures that it is uniquely
dtermined.  Therefore so is $J\!K$.

\begin{nnthm}
{\it If $\cl{Q}$ and $\tilde{\cl{Q}}$ are two quotient modules{\rm,} over the
algebra $\cl{A}(\Omega),$ in the class $B_2(\Omega, \cl{Z}),$ then they
are isomorphic if and only if}

\begin{description}
\item {\rm tan:} $\cl{K}_{\rm tan} = \tilde{\cl{K}}_{\rm tan}$
\item {\rm trans:} $\cl{K}_{\rm trans} = \tilde{\cl{K}}_{\rm trans}$
\item {\rm angle:} $\langle (\bar{\partial} e)(w) , e(w) \rangle =
\langle (\bar{\partial}\tilde{e})(w) , \tilde{e}(w) \rangle $.
\end{description}
\end{nnthm}

\begin{proof}
Suppose, we are given two quotient modules $\cl{Q}$ and
$\tilde{\cl{Q}}$ which are isomorphic.  Then the module map
$\Phi:\mathcal{Q} \to \tilde{\mathcal{Q}}$ induces an anti-holomorphic
bundle map $\Phi:J\!E_{|{\rm res}~\mathcal{Z}} \to J\!\tilde{E}_{|{\rm
res}~\mathcal{Z}}$.  For $w\in \mathcal{Z}$, let $J\!E(w)$ and
$J\!\tilde{E}(w)$ denote the two dimensional space spanned by $\{ e(w),
(\bar{\partial}e)(w)\}$ and $\{ \tilde{e}(w),
(\bar{\partial}\tilde{e})(w)\}$, respectively.  Then the bundle map
$\Phi$ defines a linear map $\Phi(w): J\!E(w) \to J\!\tilde{E}(w)$.  The
map $\Phi(w)$ must then intertwine the two nilpotents $N(w)$ and
$\tilde{N}(w)$ which implies that $\Phi(w)$ must be of the form
$\Phi(w) = \left (\begin{smallmatrix} \alpha(w) & \beta(w) \cr 0 &
\alpha(w) \end{smallmatrix}\right )$, where $\alpha, \beta$ are
anti-holomorphic functions for $w$ in some small open set in
$\mathcal{Z}$.  We observe that $\Phi(w)$ maps $\gamma_0(w)$ to
$\alpha(w) \|\tilde{e}(w)\| \|e(w)\|^{-1} \tilde{\gamma}_0(w)$.  Since
$\Phi(w)$ is an isometry, it follows that $\alpha(w)= \|e(w)\|
\|\tilde{e}(w)\|^{-1}$.  Because we have chosen to work only with
normalised kernels, we infer that $\|e(w)\| \|\tilde{e}(w)\|^{-1}=1$
for all $w\in \mathcal{Z}$ which is the same as saying that $\alpha(w)=1$ for
$w\in \mathcal{Z}$.  The condition `tan' of the theorem is evident.

The module map $\phi$ has to satisfy the relation
\begin{equation*}
  J\!K(z,w) = \overline{\Phi(z)}J\tilde{K}(z,w) \Phi(w),~~z,w \in \mathcal{Z}.
\end{equation*}
However, $J\!K(z,0)= \left (\begin{smallmatrix} 1 & 0\\ 0 & S(z)
\end{smallmatrix}\right )$, and similarly $\tilde{K}$ at $(z,0)$ has a
matrix representation with $S$
replaced by $\tilde{S}$.  Now, evaluate the formula relating $JK$ and
$J\tilde{K}$ at $w=0$ to conclude that $\overline{\beta(z)}=0$ for all
$z\in \mathcal{Z}$.

Now, since $\Phi(w)$ has to preserve the inner products, it follows
that $\langle (\bar{\partial}e)(w) , e(w) \rangle - \langle
(\bar{\partial}\tilde{e})(w), \tilde{e}(w) \rangle = \beta(w)
\|e(w)\|^2$.\ Hence it follows that $\langle
(\bar{\partial}e)(w), e(w) \rangle = \langle
(\bar{\partial}\tilde{e})(w)$, $\tilde{e}(w) \rangle$ which is the
condition `angle' of the theorem.

Finally, the requirement that the nilpotents $N(w)$ and $\tilde{N}(w)$
must be unitarily equivalent  for each $w\in \mathcal{Z}$ amounts to
the equality of the $(1,2)$ entry of $N_{\rm orth}(w)$ with that of
$\tilde{N}_{\rm orth}(w)$.  Since we have already ensured $\|e(w)\| =
\|\tilde{e}(w)\|$, it follows that $b(w) = \tilde{b}(w)$.  This  clearly
forces the condition `trans' of the theorem which completes the
proof of necessity.

For the converse, first prove that the natural map from $J\!E(w)$ to $J\!\tilde{E}(w)$,
$w\in \mathcal{Z}$, which carries one anti-holomorphic frame to the
other is an isometry.  It is evident that this map, which we denote by
$\Phi(w)$, defines an anti-holomorphic bundle map and that it
intertwines the nilpotent action.

To check if $\Phi(w)$ is isometric, all we have to do is see if it
automatically maps the orthonormal basis $\{\gamma_0(w),
\gamma_1(w)\}$ to the corresponding orthonormal basis
$\{\tilde{\gamma}_0(w), \tilde{\gamma}_1(w)\}$.  Clearly, $\Phi(w)
(\gamma_0(w)) = \tilde{e}(w)\|e(w)\|^{-1} = \tilde{\gamma}_0(w)
\|\tilde{e}(w)\|\,\|e(w)\|^{-1}$.  Suppose that the two curvatures
corresponding to the bundles $J\!E$ and $J\!\tilde{E}$ agree on the
hypersurface $\cl{Z}$. Then it is possible to find sections of these
bundles which have the same norm. Or, equivalently, we may assume that
$\|\gamma_0(w)\| = \|\tilde{\gamma}_0(w)\|$.  It then follows that
$\Phi(w)(\gamma_0(w)) = \tilde{\gamma}_0(w)$.  Notice that
\begin{align*}
\Phi(w) (\gamma_1(w)) &= a(w) \tilde{e}(w)
+ b(w) (\partial\tilde{e})(w) \\
&= a(w) \|\tilde{e}(w)\| \tilde{\gamma}_0(w) + b(w)
(\tilde{b}(w))^{-1} (\tilde{\gamma}_1(w)\\
&\quad\ - \tilde{a}(w)
\|\tilde{e}(w))\| \tilde{\gamma}_0(w)\\
&= (a(w)\tilde{b}(w) - \tilde{a}(w) b(w) )
\|\tilde{e}(w)\|(\tilde{b}(w))^{-1} \tilde{\gamma}_0(w)\\
&\quad\ + b(w) (\tilde{b}(w))^{-1} \tilde{\gamma}_1(w).
\end{align*}
A simple calculation shows that
\begin{align*}
a(w)\tilde{b}(w) - \tilde{a}(w) b(w)
&= \|e(w)\|^3 \|\tilde{e}(w)\|
(- \cl{K}(w))^{-1/2}(- \tilde{\cl{K}} (w))^{-1/2}\\
&\quad\ \big ( \langle (\bar{\partial} e)(w) , e(w) \rangle -
\langle (\bar{\partial}
\tilde{e})(w) , \tilde{e}(w) \rangle \big ).
\end{align*}
It follows that $\Phi(w)$ maps $\gamma_1(w)$ to $\tilde{\gamma}_1(w)$
if and only if $b(w) = \tilde{b}(w)$ and $\langle
(\bar{\partial} e)(w) , e(w) \rangle = \langle (\bar{\partial}
\tilde{e})(w) , \tilde{e}(w) \rangle $.

We have therefore shown that the two bundles $J\!E$ and $J\!\tilde{E}$ are
locally equivalent (via the bundle map $J\varphi$).  We now apply the
Rigidity Theorem (\cite{C-D}, p.~202) to conclude that the two modules
$\mathcal{Q}$ and $\tilde{\mathcal{Q}}$ must be isomorphic.\hfill \ab
\end{proof}

It is not clear if the condition `angle' of the theorem can be
reformulated in terms of intrinsic geometric invariants like the
second fundamental form etc.

In the case $k>2$, if we show that the bundle map is the identity
transform on each of the fibers, then it will follow that the matrix
entries of the two nilpotent actions on each of these fibers must be
equal.  These entries are expressible in terms of the curvature in the
transverse direction and its normal derivatives.  So if two quotient
modules are isomorphic, then it follows that these quantities must be
equal. However, we are not sure what a replacement for the condition
`angle' in the theorem would be.

\section*{Acknowledgements}

The second named author (GM) would like to thank
Indranil Biswas, Jean-Pierre Demailly and Vishwambhar Pati for
many hours of helpful conversations.

The research of both the authors was supported in part by DST-NSF S\&T
Cooperation Programme.  The research of the second author (GM)
was also partially funded through a grant from The National Board for
Higher Mathematics, India.

\end{document}